\documentclass[times,doublespace]{qjrms4}

\usepackage[colorlinks,bookmarksopen,bookmarksnumbered,citecolor=red,urlcolor=red]{hyperref}

\newcommand\BibTeX{{\rmfamily B\kern-.05em \textsc{i\kern-.025em b}\kern-.08em
T\kern-.1667em\lower.7ex\hbox{E}\kern-.125emX}}

\usepackage{amsmath} 
\usepackage{moreverb}
\usepackage{titlesec}

\usepackage{graphicx}

\usepackage{caption} 
\usepackage[position=top]{subfig} 
\usepackage {tikz} 
\usetikzlibrary{shapes,shadows,arrows} 
\usepackage{color} 

\begin{document}

\runningheads{A.~S.~et al.}{Assimilation of semi-qualitative observations}

\title{Assimilation of semi-qualitative observations with a stochastic Ensemble Kalman Filter}

\author{Abhishek Shah\textsuperscript{a}\corrauth, Mohamad El Gharamti\textsuperscript{a,b} and Laurent Bertino\textsuperscript{a}}

\address{\textsuperscript{a}Nansen Environmental and Remote Sensing Center, Bergen, Norway\\ \textsuperscript{b}National Center for Atmospheric Research, Colorado, USA }

\corraddr{Nansen Environmental and Remote Sensing Center, Bergen, Norway.\\ E-mail: abhishek.shah@nersc.no}

\begin{abstract}

The Ensemble Kalman filter assumes the observations to be Gaussian random variables with a pre-specified mean and variance. In practice, observations may also have detection limits, for instance when a gauge has a minimum or maximum value. In such cases most data assimilation schemes discard out-of-range values, treating them as ``not a number", at a loss of possibly useful qualitative information.

The current work focuses on the development of a data assimilation scheme that tackles observations with a detection limit. We present the Ensemble Kalman Filter Semi-Qualitative (EnKF-SQ) and test its performance against the Partial Deterministic Ensemble Kalman Filter (PDEnKF) of \citet{Borup2015}. Both are designed to explicitly assimilate the out-of-range observations: the out-of-range values are qualitative by nature (inequalities), but one can postulate a probability distribution for them and then update the ensemble members accordingly. The EnKF-SQ is tested within the framework of twin experiments, using both linear and non-linear toy models. Different sensitivity experiments are conducted to assess the influence of the ensemble size, observation detection limit and a number of observations on the performance of the filter. Our numerical results show that assimilating qualitative observations using the proposed scheme improves the overall forecast mean, making it viable for testing on more realistic applications such as sea-ice models. 

\end{abstract}

\keywords{Data Assimilation; Ensemble Kalman filter; semi-qualitative information; out-of-range observations; detection limit.}

\maketitle

\footnotetext[2]{Please ensure that you use the most up to date
class file,
available from the QJRMS Home Page at\\
\href{http://onlinelibrary.wiley.com/journal/10.1002/(ISSN)1477-870X}{\tiny\texttt{http://onlinelibrary.wiley.com/journal/10.1002/(ISSN)1477-870X}}%
}

\section{Introduction}

Data Assimilation (DA) is an approach through which available observations along with the prior knowledge (model state) are used to obtain an estimate of the true state of a process \citep{Ghil1991,Daley1993,Talagrand1997,Kalnay2003}. Each observation is used to reduce model uncertainty and improve its forecast accuracy. In practice, many observations are only available in a limited interval of the actual variation of observed quantity i.e., observations with detection limit. For instance, observations with higher detection limit are Soil Moisture and Ocean Salinity (SMOS) satellite estimates of the sea-ice thickness \citep{Kaleschke2012,Kaleschke2010} and ocean winds observations from scatterometers in hurricane wind speeds \citep{Reul2012}. SMOS can give quantitative thickness data only up to 50cm over first-year level ice for the Arctic, because the signal penetration is limited by the wavelength. In reality, the sea ice can grow up to a few meters. Conversely, observations with lower detection limit also exist. Examples are contaminant concentrations in environmental and health fields \citep{Hornung1990} and river water level measurements obtained from satellite radar altimetry. On top of detection limits, some measurements are boolean in nature for e.g. if the permafrost exists or not \citep{Li1999}, or whether or not there is overflow at a weir in urban hydrology \citep{Thorndahl2008}. Although these types of observations do not provide quantifiable data above or below the detection limit, they do give qualitative information about the observed variable. Therefore, this type of observations should be exploited as a means to improve the model forecast. 

All deterministic and stochastic ensemble based filtering schemes \citep{Burgers1998,Anderson2001,Tippett2003,Sakov2008} assimilate actual observations (hard data), but do not consider qualitative information (soft data) available from out-of-range observations (OR-observations). Whereas the geostatistical techniques are well established for variables without dynamical evolution \citep{Chiles1999,Emery2008}, only one study \citep{Borup2015}, to the best of our knowledge, has dealt with the issue of OR-Observations in a dynamical data assimilation framework. 

\citet{Borup2015} proposed the Partial Deterministic Ensemble Kalman Filter (PDEnKF) to assimilate observations with a detection limit. The main idea of the PDEnKF is to assume a virtual observation at the detection limit in the absence of hard data and defining a constant OR-observation likelihood in unobservable region from the detection limit. The virtual observation is then used to update the anomalies within the framework of the Deterministic Ensemble Kalman Filter (DEnKF) \citep{Sakov2008}. Anomalies are updated differently conditioned on the values of forecast ensemble members i.e., whether the member is inside or outside the observable range. The mean, on the other hand, is updated only when there is a hard data or else there is no update. Practically, the virtual OR-observation is used only to update the ensemble members that are within the observable range. The scheme has been tested using both linear and non-linear reservoir cascade models. The authors present important improvement in the forecasts accuracy, implying that soft data can contribute meaningful information to predictions. 

In light of this background, a new DA algorithm referred to as Ensemble Kalman Filter Semi-Qualitative (EnKF-SQ) is developed here and is designed to explicitly assimilate the OR-observations. \textcolor{black}{The EnKF-SQ assumes a virtual observation at the detection limit in the absence of hard data, with an asymmetric two-piece Gaussian observational
likelihood on either side of the detection limit.} In the EnKF-SQ, the forecast ensemble members are updated by the observations, which are perturbed using a two-piece Gaussian OR-observation likelihood following the stochastic Ensemble Kalman Filter (EnKF) update \citep{Evensen2003}. Detailed derivation of EnKF-SQ is discussed in Section \ref{subsec:EnKF-SQ} followed by an algorithmic implementation. To test the performance of the EnKF-SQ, we apply it to two different linear and non-linear toy models. The experimental setup and results are presented in section \ref{sec:Numerical-tests}. A summary of the numerical results is followed by a general discussion that concludes the paper in section \ref{sec:Conclusion}.

\section{Methodology and Algorithm\label{sec:Methodology-and-Algorithm}}

In this section, a brief background on the stochastic EnKF is given. The new EnKF-SQ is further derived and presented in details.

\subsection{Background\label{subsec:Background}}

The Kalman Filter (KF) \citep{Kalman1960} is a sequential filtering technique, in which the model is integrated forward in time and, when they become available, observations are used to update the model state and its associated uncertainty. The KF is a recursive Bayesian estimation method, which is optimal for Gaussian and linear models \citep{Gharamti2011}. The KF operates sequentially in time following time update (forecast) and measurement update (analysis) steps. The EnKF, a variant of the KF, utilizes an ensemble of model states $[\mathbf{x}_{1},\mathbf{x}_{2}.....\mathbf{x}_{N}]$ (where $N$ is the ensemble size) to estimate the mean and covariance. The analysis step of the EnKF at any particular time is given as 
\begin{equation}
\mathbf{x}_{i}^{a}\mathbf{=x}_{i}^{f}\mathbf{+K}\mathbf{(}\mathbf{y}_{i}\mathbf{-H}\mathbf{x}_{i}^{f}\mathbf{)},\hspace{1.5em}i=1,2,....,N,\label{eq:EnKF_update_eq} 
\end{equation} 
\begin{equation} 
\mathbf{K}=\mathbf{P}^{f}\mathbf{H}^{T}(\mathbf{H}\mathbf{P}^{f}\mathbf{H}^{T}+\mathbf{R})^{-1},\label{eq:Kalman_gain} 
\end{equation}
where $\mathbf{K}$ is referred to as the Kalman gain; $\mathbf{x}_{i}^{a}$ and $\mathbf{x}_{i}^{f}$ is the $i^{th}$ analysis and forecast state member, respectively; $\mathbf{y}_{i}$ is the $i^{th}$ vector of perturbed observations; $\mathbf{H}$ is the observation operator i.e., mapping the state variable to the observation space (assumed linear here); $\mathbf{P}^{f}$ is the ensemble forecast error covariance matrix and $\mathbf{R}$ is the observation error covariance matrix. The superscripts $a$, $f$ and $T$ stand for analysis, forecast and matrix transpose, respectively. For clarity, the time index is omitted from all notations. The term $\mathbf{(}\mathbf{y}_{i}\mathbf{-H}\mathbf{x}_{i}^{f}\mathbf{)}$ in eq.~(\ref{eq:EnKF_update_eq}) is the discrepancy between the observations and the ensemble members, often referred to as sample innovations. The ensemble forecast error covariance matrix $\mathbf{P}^{f}$ is never explicitly computed, however, it is decomposed as follows: 
\begin{equation} 
\mathbf{P}^{f}=\frac{1}{N-1}\sum\limits_{i=1}^N(\mathbf{x}_{i}^{f}-\bar{\mathbf{x}})(\mathbf{x}_{i}^{f}-\bar{\mathbf{x}})^{T}=\frac{1}{N-1}\mathbf{AA}^{T}, 
\end{equation}
where $\bar{\mathbf{x}}$ is the mean of the forecast ensemble and $\mathbf{A}=[\mathbf{A}_{1},\mathbf{A}_{2}......\mathbf{A}_{N}]$ is the ensemble anomalies matrix. Similarly, the analysis error covariance matrix can be computed from the ensemble of analysis states, but is not required in the implementation.

For hard data with known observational likelihood, each ensemble member is updated independently using observations that are perturbed with ${\cal N}(0,\mathbf{R})$ as shown in \citet{Burgers1998,Evensen2003}. For observations with detection limits, the likelihood is truncated and is therefore non-Gaussian. How can this information be incorporated in an EnKF system?

\subsubsection{Observations with a detection limit \label{subsec:Observation-with-detection}}

When the observations have a detection limit, one may not have a full access to the observation likelihood. For simplicity, we will only consider the case with an upper detection limit on the observations rather than a lower limit, without loss of generality. Observations with detection limit can be characterized in two parts: 
\begin{enumerate}
\item Hard data or in-range observations $(\mathbf{y}_{\mathit{ir}})$. 
\item Soft data or OR-observations $(\mathbf{y}_{\mathit{or}})$ i.e., no
specific value of the observed quantity.
\end{enumerate}
Within the Bayesian framework, the goal of DA is to estimate the posterior distribution of the model state. According to Bayes' rule, the posterior distribution $p(\mathbf{x|y})$ is proportional to the product of a prior $p(\mathbf{x})$ and the observation likelihood $p(\mathbf{y|x})$ as follows: 
\begin{equation}
p(\mathbf{x|y})\propto p(\mathbf{y|x})p(\mathbf{x})\label{eq:Bayes_rule}
\end{equation}

For an observation with detection limit, eq.~(\ref{eq:Bayes_rule})
can be split into two, depending on the nature of the observation,
i.e., 
\begin{equation}
p(\mathbf{x|y}) \propto
\begin{cases}
p(\mathbf{y_{\mathit{ir}}}|\mathbf{x})p(\mathbf{x}),\hspace{0.5cm}\text{when \ensuremath{\mathbf{y=y}_{ir}} in-range observations}\\
p(\mathbf{y_{\mathit{or}}}|\mathbf{x})p(\mathbf{x}),\hspace{0.45cm}\text{when \ensuremath{\mathbf{y=y}_{or}} OR-observations}
\end{cases}\label{eq:Bayes_yir_yor}
\end{equation}

Although we do not have a priori OR-observation likelihood $p(\mathbf{y}_{or}|\mathbf{x})$ to solve for the posterior $p(\mathbf{x}|\mathbf{y})$ given in eq.~(\ref{eq:Bayes_yir_yor}), one can always postulate a OR-observation likelihood based on climatology or experts opinion. A detailed discussion about the choice of OR-observation likelihood is given in section \ref{subsec:EnKF-SQ}. In the following section we will introduce the EnKF-SQ and present its implementation along with its main differences to the PDEnKF. 

\subsection{The Ensemble Kalman Filter Semi-Qualitative\label{subsec:EnKF-SQ}}

The KF minimizes the forecast error variance, and this is achieved by updating state variables, eventually moving them closer to the observations. The update of the prior given in eq.~(\ref{eq:Bayes_yir_yor}) for in-range observations is straightforward. But for OR-observations it is not so clear, since we do not have distribution of the observation. As such, an assumption about the OR-observation likelihood, which should be physically consistent with the observed quantity and the qualitative information we have about it, is required to solve eq.~(\ref{eq:Bayes_yir_yor}). 

\subsubsection{The Partial Deterministic EnKF}

As a way to do that, \citet{Borup2015} proposed a DA scheme, namely PDEnKF, to solve the Bayesian system in eq.~(\ref{eq:Bayes_yir_yor}). The authors assumed the OR-observation likelihood to be constant outside the observable range. Furthermore, the likelihood function inside the observable range is set to be determined by in-range observation uncertainty because measurement errors make it possible for in-range values to be wrongly observed as out of range. Inspired by the PDEnKF, we present the EnKF-SQ that uses a stochastic EnKF.

In contrast to the stochastic update, the PDEnKF follows \citet{Sakov2008} and uses two different equations for updating the ensemble mean and anomalies, the anomalies being updated by half the gain in a form of implicit inflation. In some cases of partial update, the half-gain does not maintain the anomalies centered on the analysis mean: if the mean is within the range and the observation outside, the ``half-gain'' will leave the anomalies further inside the range than if the partial update were applied to the mean. \citet{Borup2015} have opted for this non-centered partial analysis scheme in order to maintain more ensemble spread.

In the EnKF-SQ, instead of a constant uniform OR-observation likelihood, we propose to use a two-piece Gaussian distribution \citep{Gibbons_two_piece_1973,Fechner_2piece} as OR-observation likelihood. In other words, the uniform likelihood of \citet{Borup2015} is replaced by a Gaussian distribution with varying observation error variance outside the observable range. 

A two-piece Gaussian distribution is obtained by merging two opposite halves of the two Gaussian probability densities (pdfs) at their common mode, given as follows: 
\begin{equation}
f(x)=\begin{cases}
W\mathbf{exp}\left[-\frac{\left(x-\mu\right)^{2}}{2\sigma_{1}^{2}}\right],\hspace{0.5cm}x{\leq}\mu\\
W\mathbf{exp}\left[-\frac{\left(x-\mu\right)^{2}}{2\sigma_{2}^{2}}\right],\hspace{0.5cm}x{>}\mu
\end{cases}
\end{equation}
where $W\mathbf{=\sqrt{\frac{2}{\pi}}\left(\mathbf{\sigma_{1}+\sigma_{2}}\right)^{-1}}$ is a normalizing constant, $\mathbf{\mu}$ is the common mean, $\mathbf{\sigma_{1}}$ and $\sigma_{2}$ are the standard deviations (std) of the two Gaussian pdfs. The common mean $\mu$ is located at the detection limit, as it is the last possible value the gauge could detect with known observation uncertainty. In essence, the common mean $\mu$ is nothing but the mode of a two-piece Gaussian distribution. Note that for the function $f(x)$, the mean does not coincide with the mode given the skewness of the distribution. The reasons for choosing a Gaussian likelihood, over a uniform one, in the unobservable range are:

\begin{itemize}
\item OR-observations do not give a specific value of the observed quantity, but an educated guess can always be made about a realistic range of values. For instance using a climatology of the values in the unobservable range. Imposing a uniform density outside the observable range gives equal weight to all values until infinity, whereas extremely high values are usually less realistic in most applications like wind speed, ice thickness among others.
\item In order to implement the stochastic EnKF, one needs to perturb the observations eq.~(\ref{eq:EnKF_update_eq}) with a Gaussian distribution of covariance matrix $\mathbf{R}$. For the OR-uniform likelihood it is technically impossible since the uniform tail is not integrable. Even if the OR uniform likelihood were limited to a finite upper bound, the choice of that upper bound would have to be justified by the nature of the variable. In principle there is no restriction to the choice of OR likelihood probability distribution but the Gaussian distribution has practically convenient properties for our purpose (the simulation do not generate excessive outliers and the Bayesian interpretation is relatively simple, see below).
\end{itemize}

\begin{figure}
\begin{centering}
\includegraphics[width=14cm,height=8cm]{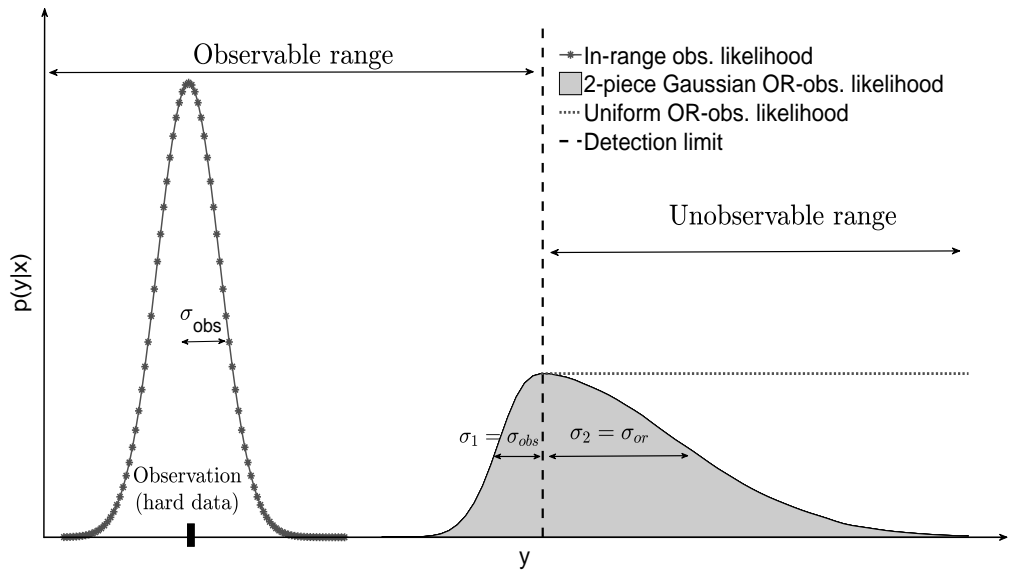}
\par\end{centering}
\begin{centering}
\caption{Illustration of the two-piece Gaussian OR-observation likelihood, when a gauge has an upper observation limit. In-range observation is shown in small black rectangle and the corresponding Gaussian likelihood in solid black line. The two-piece Gaussian likelihood is plotted in grey and Uniform OR-obs. likelihood from \citet{Borup2015} is in dashed-grey. $\sigma_{obs}$ is the observation error std for hard data and $\sigma_{or}$ is the educated guess of the OR-obs error std.\label{fig:Illustration_2_piece_Gaussian_likelihood} }
\end{centering}
\end{figure}

In addition, we assume that the observation error variance of the Gaussian half which is inside the observable range from detection limit is equal to the in-range observation error standard deviation $(\sigma_{obs})$ as in \citet{Borup2015}. An example of the two-piece Gaussian OR-observation likelihood having an upper detection limit is shown in Fig.~\ref{fig:Illustration_2_piece_Gaussian_likelihood}. As shown, the two-piece Gaussian likelihood is right-skewed because of the higher OR-observation error standard deviation $\sigma_{or}$. Choosing a \textit{proper} $\sigma_{or}$ is very important as it will be used to generate perturbations and thereby to update ensemble members. The choice should be in adequacy with the possible values in the unobservable range of the underlying observed variable. 

\subsubsection{Choice of $\mathbf{\sigma}_{or}$\label{subsec:Choice-of}}

The observation error standard deviation for the Gaussian half outside the observable range $(\mathbf{\sigma}_{or})$ is an arbitrary choice with different possibilities. If the pdf of the climatological data of the observed quantity is available, then $\sigma_{or}$ can be approximated by using the mean of out-of-range climatological values:
\begin{equation}
\mathbf{\sigma}_{or}=-\mu+\left(\int\displaylimits_{\mu}^{+\infty}\mathbf{y}f_{clim}(\mathbf{y})dy\right)\label{eq:sigma_OR}
\end{equation}
where $f_{clim}(\mathbf{y})$ is the pdf of the climatological data of the observed quantity, $\mu$ is the detection limit point. The second term on the right hand side of eq.~(\ref{eq:sigma_OR}) is the expectation of the climatological distribution for the values above the detection limit. Eq.~(\ref{eq:sigma_OR}) is used to generate $\sigma_{or}$ values in all of the experiments presented in section \ref{sec:Numerical-tests}. Sensitivity experiments using different values for $\sigma_{or}$ are also conducted (section \ref{subsec:Sensitivity_experiments}). In the absence of climatology for the observed data, an educated guess can be used based on expert knowledge about the $\sigma_{or}$, considering that extremely high values are less likely and vice versa for lower detection limit. 

\subsubsection{Bayesian representation \label{subsec:Bayesian-representation}}

According to Bayes' rule, the posterior distribution is proportional to the product of the prior and observation likelihood functions eq.~\ref{eq:Bayes_rule}. For hard data, the posterior is simply the product of two Gaussian distributions and it is Gaussian. For OR-observations, it is the product of a Gaussian prior distribution and a two-piece Gaussian likelihood. This is nothing but the product of two Gaussian distributions (\textit{the prior and each half of a two-piece Gaussian}) on either side of the detection limit and hence the posterior is a piecewise Gaussian distribution meeting at $\mu$. By construction, the posterior is unimodal.

Fig.~\ref{fig:bayesian_representation_illustatiom} illustrates the update when (i) the mode of the prior distribution is inside and (ii) the mode of the prior distribution is outside the observable range. The curve for two-piece Gaussian likelihood and posterior-Bayes' distribution in Fig.~\ref{fig:bayesian_representation_illustatiom} are obtained by sampling the respective pdf with the acceptance-rejection method \citet{VonNeumann1951} of Monte Carlo techniques. We have assumed a Gaussian distribution ${\cal N}(\mu,2\sigma_{or}^{2})$ as a proposal distribution to generate the samples using the acceptance-rejection method. We have generated 100000 samples to plot the figure. 

As shown, when the mode of a prior distribution is inside the observable range (Fig.~\ref{fig:update_mem_inside}), the location of the posterior mode will be between the mode of the prior and the observation detection limit. This demonstrates the desired effect on the prior distribution by moving it towards the unobservable range. If the mode of the prior distribution is outside the observable range, then the update has a very small effect on it (Fig.~\ref{fig:update_mem_outside}), as expected because of the high observation error variance outside of the observable range. 

\subsubsection{Implementation and Algorithm} \label{subsec: algorithmic_implementatin}

The ensemble members can be seen as discrete samples of a continuous distribution. Updating the ensemble members given the hard data is performed by the stochastic EnKF as in eq.~(\ref{eq:EnKF_update_eq}). As for the soft data, the Bayes' equation shows that ensemble members inside the observable range needs to be updated towards the unobservable range. Intuitively, the ensemble members which lie in the unobservable range should be left untouched, as we do not have a specific value of the observation. 

The proposed approach to update the ensemble for soft data is divided into two cases whether the observed ensemble members i.e.~$\mathbf{Hx}_{i}^{f}$ are inside or outside the observable range. Members inside (outside) the observable range should be updated linearly with observation uncertainty $\sigma_{obs}$ ($\sigma_{or}$). During the update, the observation perturbations $\mathbf{y}_{i}$ can be generated by the acceptance-rejection method assuming a Gaussian function as a proposal distribution with mean at $\mu$(detection limit) and standard deviation equal to the maximum of $\sigma_{obs}$ and $\sigma_{or}$; i.e. $max(\sigma_{or},\sigma_{obs})$, so that it has longer tails than the two-piece Gaussian likelihood. 
\begin{figure}
\subfloat[Mode of a prior is inside the range\label{fig:update_mem_inside} ]{\includegraphics[width=8.3cm,height=6.2cm]{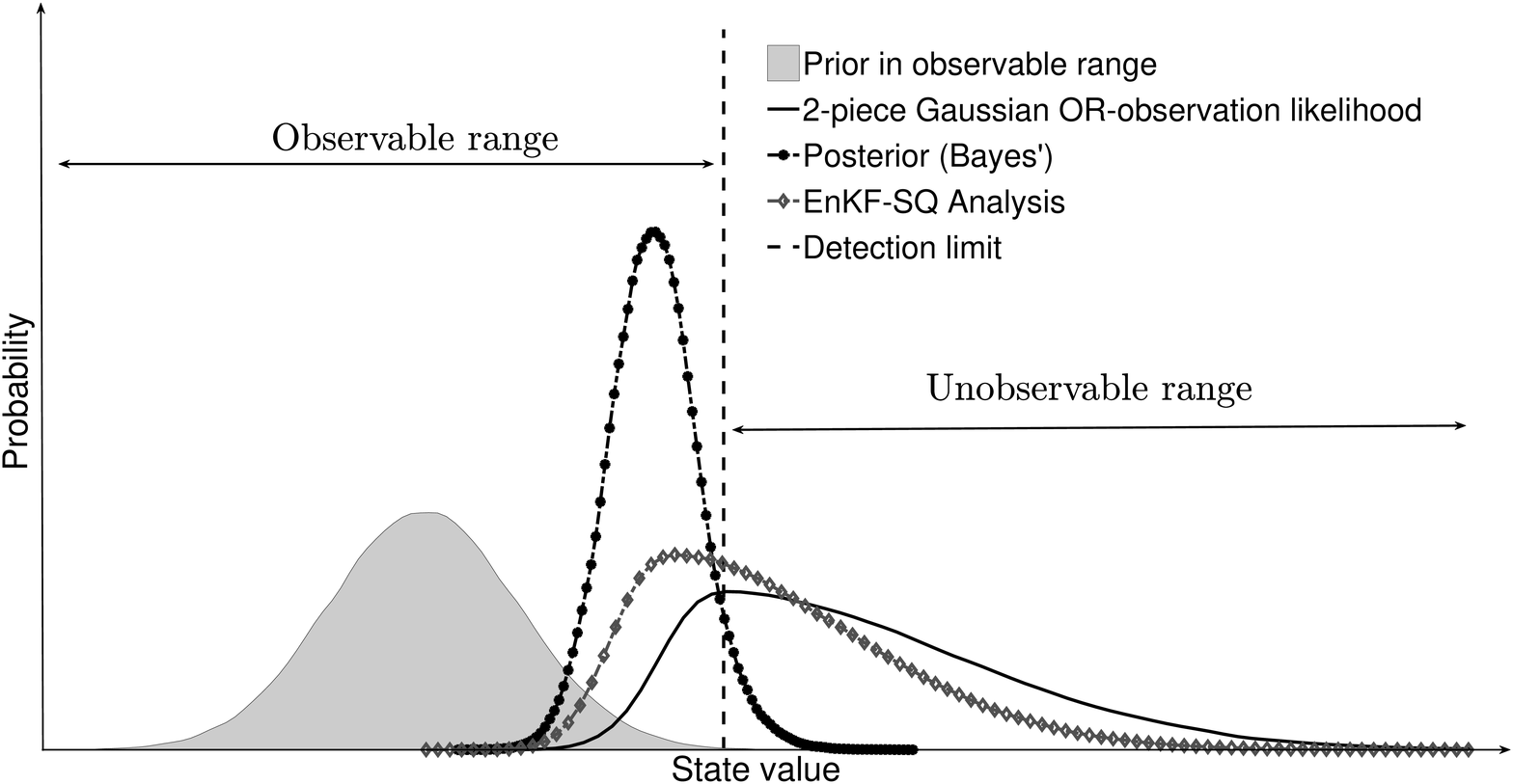}}
\hspace{1em}
\subfloat[Mode of a prior is outside the range\label{fig:update_mem_outside} ]{\includegraphics[width=8.3cm,height=6.2cm]{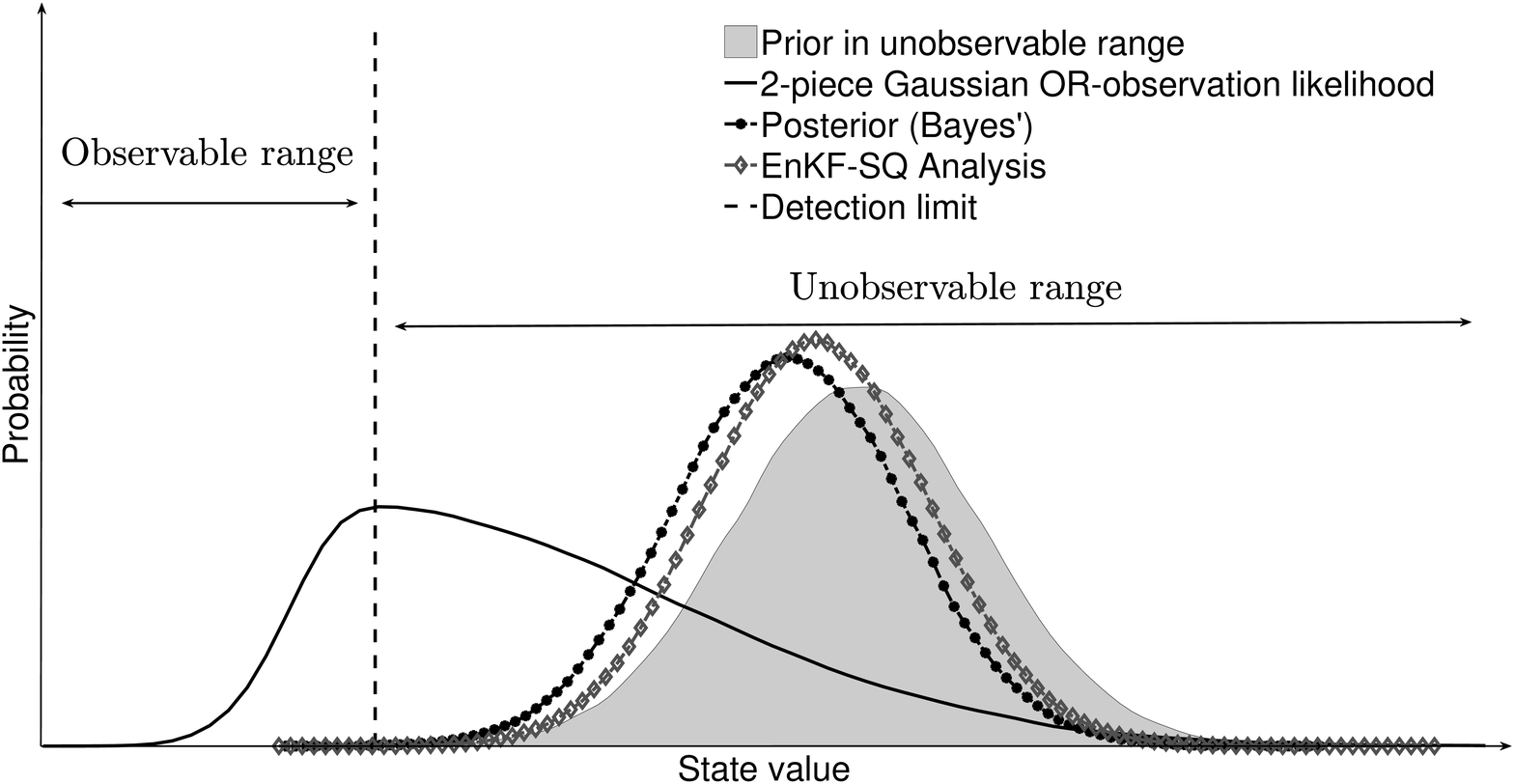}}\caption{Bayesian posterior and the EnKF-SQ analysis for a scalar update of a Gaussian prior with a two-piece Gaussian OR-observation likelihood. (a) the mode of a prior is inside the observable range. (b) the mode of a prior is in the unobservable range. \label{fig:bayesian_representation_illustatiom}}
\end{figure}
The Kalman gain $\mathbf{K}$ for the forecast ensemble member inside the observable range is calculated with in-range observation error standard deviation $\sigma_{obs}$. For a scalar case, this becomes $K_{ir}=\sigma_{b}^{2}\left(\sigma_{b}^{2}+\sigma_{obs}^{2}\right)^{-1}$. If the member is outside the observable range then the Kalman gain $\mathbf{K}$ is calculated with out-of-range observation error variance $\mathbf{\sigma}_{or}$ i.e., $K_{or}=\sigma_{b}^{2}\left(\sigma_{b}^{2}+\sigma_{or}^{2}\right)^{-1}$. Here, the forecast error variance is denoted by $\sigma_{b}^{2}$. 

For multivariate case this can be achieved by simply changing the values of observation error variance to $\sigma_{obs}^{2}$ or $\sigma_{or}^{2}$ for OR-observation in the error covariance matrix $\mathbf{R}$, depending on the location of forecast ensemble member. Note that the proposed algorithm for EnKF-SQ only supports uncorrelated observations i.e.~matrix
$\mathbf{R}$ is diagonal. If the observations errors are correlated, one can decorrelate them \citep{Evensen2004} and proceed with the algorithm. An algorithmic implementation of the EnKF-SQ analysis is presented below:

\begin{description}
\item [{Algorithmic steps:}]
\end{description}
For an efficient processing of the update eq.(\ref{eq:Bayes_yir_yor}), observations are pre-processed serially to sort out hard data ($y_{ir}$) and soft data ($y_{or}$) before proceeding to analysis. The subscripts $ir$ and $or$ stands for the index number of any hard and soft data in observation vector $\mathbf{y}$, respectively. 

\noindent For each ensemble member $i$:
\begin{enumerate}

\item For each OR-observations $y_{or}$, apply the $or^{th}$ observation operator row $\mathbf{H}_{or}$ to ensemble member $\mathbf{x}_{i}^{f}$, to check whether the member is outside or inside the observable range. 

\item Perform the operation below for all OR-observations, in order to set the values of observation error variance in matrix $\mathbf{R}$ depending on the location of $\mathbf{H}_{or}\mathbf{x}_{i}^{f}$.\\
Pseudo-code:\\
for each OR-observations $y_{or}$\\
\phantom{xx}\textit{if} $\mathbf{H}_{or}\mathbf{x}_{i}^{f}>\mu$\\
\phantom{xx xx}$\mathbf{R}_{(or,or)}=\sigma_{or}^{2}$\\
\phantom{xx}\textit{else}\\
\phantom{xx xx}$\mathbf{R}_{(or,or)}=\sigma_{obs}^{2}$\\
\phantom{xx}\textit{end if}\\
end for each $\mathbf{y}_{or}$
\item Calculate the Kalman gain matrix $\mathbf{K}$ with the updated $\mathbf{R}$.
\item Update the forecast ensemble member $\mathbf{x}_{i}^{f}$ using EnKF update eq.~(\ref{eq:EnKF_update_eq}), where the perturbation vector $\mathbf{y}_{i,or}$ and $\mathbf{y}_{i,ir}$ are generated from two-piece Gaussian likelihood and ${\cal N}(y_{ir},\sigma_{obs}^{2})$ for OR and in-range observations respectively. 
\item Repeat the process for all $N$ ensemble member $\mathbf{x}_{i}^{f}$ to get the analysis ensemble.
\end{enumerate}
End the loop on $i$.

\begin{figure}

\tikzstyle{decision} = [diamond, draw=black, thick, fill=gray!3,text width=4.5em, text badly centered, node distance=3cm, inner sep=0pt,minimum height=3em,line width =1.2]   
\tikzstyle{block} = [rectangle, draw=black, thick, fill=gray!3, text width=5em, text centered, rounded corners, minimum height=3.7em,line width = 1.2]   
\tikzstyle{line} = [draw, -latex',thick,line width = 1.2]   
\tikzstyle{cloud} = [draw=black,thick, ellipse,fill=gray!3, text centered, node distance=3cm, minimum height=2em,line width = 1.2]

\large  
\begin{center}  
\begin{tikzpicture}[node distance = 3cm, auto]    
\node [cloud,fill=gray!3,minimum height=4em, text width=7em] (start) {Update Scheme of EnKF-SQ}; 
\node [block, below of=start, text width=13em, minimum height=3em] (pre-process) {Pre-process all obs. to separate hard and soft data. \\ for each $\mathbf{x}_i^f$}; 
\node [decision, below of= pre-process,node distance=4cm, text width=5em] (check-ens) {$\textbf{H}_{or}\textbf{x}_i^f$ is within obs. range}; 
\node [block, right of=check-ens,node distance=5cm,text width=8em] (out-range-ens) {Set the value of $\textbf{R}_{(or,or)}$ to $\sigma_{or}^{2}$}; 
\node [block, below of=check-ens, node distance=3.5cm, text width=8em] (in-range-ens) {Set the value of $\textbf{R}_{(or,or)}$ to $\sigma_{obs}^{2}$}; 
\node [block, below of=in-range-ens, node distance=2.8cm, text width=13em] (perturb soft) {Perturb $y_{or}$ with two-piece Gaussian likelihood.\\Follow step 2 of algorithm for each OR-obs. $y_{or}$}; 
\node [block, below of=perturb soft, node distance=2.8cm,text width=13em,minimum height=3em] (repeat-1) { Follow the step numbers 3, 4 and 5 of algorithm.\\ End}; 
\node [block, left of=repeat-1, above=6cm, text width=8em, node distance=5.5cm] (hard data) {Perturb $y_{ir}$\\ with $\cal{N}$($y_{ir},\sigma_{obs}^{2}$)}; 

\path [line,thick] (start) -- (pre-process); 
\path [line,thick] (pre-process) -- node {\textit{Soft data}} (check-ens);  
\path [line,thick] (pre-process) -| node[near start,above] {\textit{Hard data}} (hard data) ;  
\path [line,thick] (check-ens)  -- node {\textit{No}} (out-range-ens);  
\path [line,thick]  (check-ens) -- node{\textit{Yes}}  (in-range-ens);  
\path[line,thick] (in-range-ens) -- (perturb soft); 
\path[line,thick] (out-range-ens) |- (perturb soft); 
\path[line,thick] (perturb soft) -- (repeat-1); 
\path[line,thick] (hard data) |- (repeat-1);  
\end{tikzpicture}  
\end{center} 

\caption{Flowchart for the implementation of the EnKF-SQ. Note that the algorithm does not lend itself to matrix multiplications as in \citet{Evensen2003}\label{fig:Flowchart-EnKF_SQ} }

\end{figure}
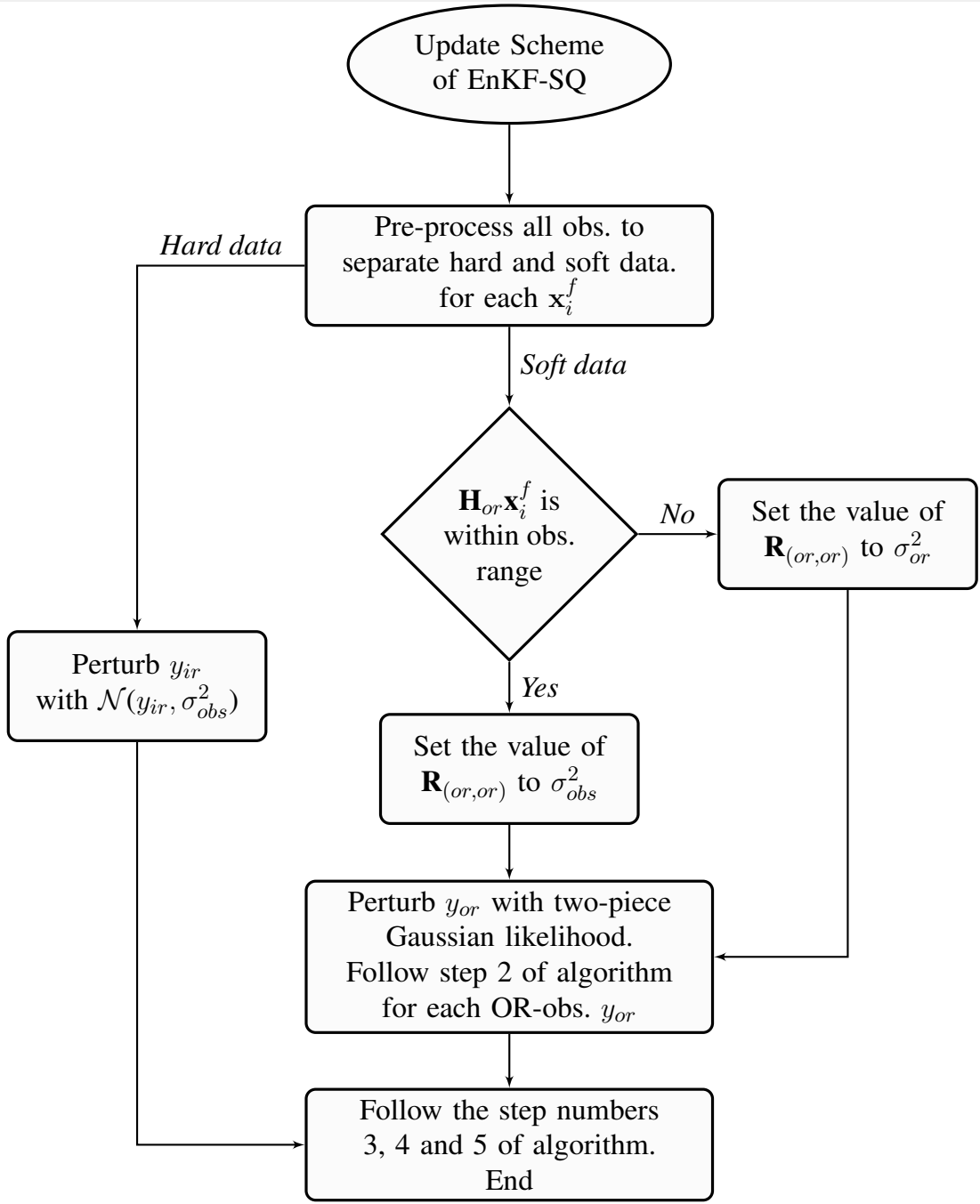

A flowchart for the EnKF-SQ update scheme is given in Fig. \ref{fig:Flowchart-EnKF_SQ}. 

To study the posterior obtained by the proposed EnKF-SQ algorithm, we superimpose the EnKF-SQ analysis to the Bayesian solution in both panels of the Fig.~\ref{fig:bayesian_representation_illustatiom}. The EnKF-SQ analysis is obtained from the exact same prior and likelihood, which as in section \ref{subsec:Bayesian-representation}. Since the likelihood is not Gaussian, we do not expect the EnKF-SQ ensemble to coincide with the Bayesian solution. We used 10000 ensemble members to sample the prior and two-piece Gaussian OR-observations likelihood, note that a more efficient algorithm then the classical acceptance-rejection may be used for computational efficiency. 

When the mode of a prior is outside of the observable range the EnKF-SQ scheme yields approximately the same posterior as that of Bayes' rule (Fig.~\ref{fig:update_mem_outside}) only marginally closer to the prior distribution. The EnKF-SQ slightly \emph{under-assimilate} in this case, which conforms to the intention of little impact of OR-observations on OR forecast members.

On the contrary, when the mode of the prior is inside the observable range the posterior obtained by Bayes' rule has a sharper peak whereas the analysis obtained from the EnKF-SQ has a thick tail in the OR domain (Fig.~\ref{fig:update_mem_inside}). The large deviation from the Bayesian solution is a sign of the sensitivity of the linear EnKF-SQ update to a skew input likelihood. This skewness was already present in the previous case, however, it was not visible with larger OR-observation errors. In the present case, the posterior EnKF-SQ ensemble is closer to the likelihood than the Bayesian solution, so it can be stated that the EnKF-SQ \emph{over-assimilate} in this case although it does return a larger ensemble spread than the Bayesian solution, which may be counter-intuitive for EnKF practitioners. It is worth noticing that the posterior modes obtained from Bayes' rule and EnKF-SQ analysis still remain close to each other as intended in \citet{Borup2015}.

Note that the posterior represented in Fig.~\ref{fig:update_mem_inside} and \ref{fig:update_mem_outside} may not be very well sampled in practice if the ensemble is very small. The inconvenience of sampling errors and skewness will be evaluated with toy models in the following sections.


\section{Numerical tests\label{sec:Numerical-tests}}
In this section we present and analyze the assimilation results obtained using the proposed EnKF-SQ algorithm. We use two different toy models to test and evaluate the behavior of the EnKF-SQ. The first is a linear subsurface flow model (LSST) and the second is the non-linear Lorenz-40 (L40) model of \citet{Lorenz1998}. We conduct various sensitivity experiment with variable ensemble size, detection limit and $\sigma_{or}$. We also compare the performance of the EnKF-SQ against the PDEnKF and with two different versions of the stochastic EnKF, denoted as follows:
\begin{enumerate}
\item EnKF-ALL: No observation detection limit is applied during DA experiments, thus all observations are hard data.
\item EnKF-IG: Assimilating only hard data and ignoring soft data during the analysis. The goal for testing with EnKF-IG is to assess the added information introduced by EnKF-SQ.
\end{enumerate}
First we give a brief description of the models and the configuration used in the tests, and then we discuss the results from different numerical experiments. 

\subsection{The Linear Subsurface Transport (LSST) Model }

We consider 1D subsurface transport model in an unconfined aquifer. The transport model is driven by a steady subsurface flow using a combined Darcy's law and Continuity equation. Groundwater flows from west to east at a reference Darcy velocity of $1.18\times 10^{-4}$ m/s. Periodic water head boundary conditions are assumed. The domain is uniformly discretized into 100 cells with each cells measuring 10 m in length. 

The generalized 1D linear solute transport model is obtained from the mass conservation of species defined as: \begin{equation} r_c\frac{\partial \left( \phi \mathbf{C} \right)}{ \partial t} + \frac{\partial \left( \mathbf{UC}\right)}{dx} = \mathbf{q}, \end{equation} where $r_c$ is the retardation coefficient, $\phi$ is the porosity, $t$ is time (s), $\mathbf{C}$ is the concentration of contaminant species (ppm) and $\mathbf{q}$ is the contaminant source. An initial condition for the concentration is specified $\mathbf{C}\left( x, 0 \right) =3+\sin(5x_i)$ where $x_i$ is the length of the $i^{th}$ grid cell. The time step is set to 10 hours. The porosity is uniform and equal to $33.4$\% with a retardation coefficient of $5.19$. The contaminant source, $\mathbf{q}$ at every time step is equal to $3\times 10^{-6}$ ppm. Water flowing from western boundary is contaminated with $5$ ppm concentration value. Using these parameters, a reference run solution is simulated for a period of 4 years. In order to mimic realistic scenarios, we impose model error in the forecast model. Essentially, we perturb the transport parameters such that $\phi = 30$\% and $r_c = 6.87$. We further add a Gaussian noise $\cal{N}$(0, 0.01) to the contaminant source and the Darcy velocity field. 

\subsection{The L40 model}

The L40 model \citep{Lorenz1998} is a chaotic and non-linear model with $40$ state variables. It imitates the evolution of an unspecified scalar meteorological quantity for instance temperature or vorticity along a latitude circle. This model has been used for testing ensemble based assimilation methods in a number of earlier studies e.g., \citep{Anderson2001,Whitaker2002,Sakov2008}. The model assumes cyclic boundary conditions as follows: 
\begin{equation}
\frac{dz_{i}}{dt}=\left(z_{i+1}-z_{i-2}\right)z_{i-1}-z_{i}+F,\quad i=1,......,40;
\end{equation}
\begin{equation*}
z_{0}=z_{40},\quad z_{-1}=z_{39},\quad z_{41}=z_{1}.
\end{equation*}
where $\mathbf{z}_{i}$ is the $i^{th}$ state variable and $F$ is a forcing term.The time step is set to $\Delta t=0.05$ units (i.e., 6 hours in real atmospheric time). The model is integrated forward in time using a fourth-order Runge-Kutta. The reference (truth) trajectory is initialized by setting $F=8$, $z_{i}=F$ $\forall$$i\neq20$ and $z_{20}=F+0.001$. A reference run solution is simulated for a period of $5$ years ($7300$ steps). Initial ensemble members are obtained by perturbing the mean state of the reference trajectory with a white noise equal to ${\cal N}(0,3)$. Observations are collected from the reference trajectory and then contaminated using ${\cal N}(0,1)$. We impose a model error for data assimilation experiments by changing the forcing parameter to $F=8.1$.

\subsection{Results\label{subsec:Results}}

Experiments are performed over a period of $5$ and $4$ years for the L40 and the LSST model, respectively. The size of the ensemble is chosen based on a series of sensitivity experiments and set to $75$ and $30$ for the L40 and LSST models, respectively. The choice is made so that tunning parameters such as inflation and localization are not needed. The goal is to assess the performance of the EnKF-SQ, PDEnKF and EnKF-IG schemes for a large enough ensemble without the necessity to mitigate sampling errors and other filter related deficiencies. For the L40 model, all $40$ variables are observed and assimilated every day (i.e., every fourth time step). In the LSST model, $80$ variables are observed, with a regularly spaced observing network, every tenth time step. 

The forecast root mean square error (RMSE) is used to evaluate the filters' performances. Given the $n-$dimensional mean forecast state vector $\hat{\mathbf{x}}_{t}^{f}=(\hat{x}_{1,t}^{f},\hat{x}_{2,t}^{f},...,\hat{x}_{n,t}^{f})$ at time $t$ and if $t_{max}$ is the final time, then the time-averaged RMSE is defined as: 
\begin{equation}
\widehat{RMSE}=\frac{1}{t_{max}}\sum\limits_{t=1}^{t_{max}}\sqrt{\frac{1}{n}\sum\limits_{i=1}^n(\hat{x}_{t,i}^{f}-x_{t,i}^{r})^{2}}
\end{equation}
where $\mathbf{x}_{t}^{r}=(\hat{x}_{1,t}^{r},\hat{x}_{2,t}^{r},...,\hat{x}_{n,t}^{r})$ is the reference state vector at time $t$. Each filter run is then repeated $L=10$ times, with different random seeds to initialise the random number generator. The average RMSE over these $L$ runs is then reported as: 
\begin{equation}
\overline{RMSE}=\frac{1}{L}\sum\limits_{l=1}^L\widehat{RMSE_{l}}
\end{equation}


\subsubsection{General behavior of the EnKF-SQ}

\begin{figure}
\begin{centering}
\includegraphics[width=17cm,height=15cm]{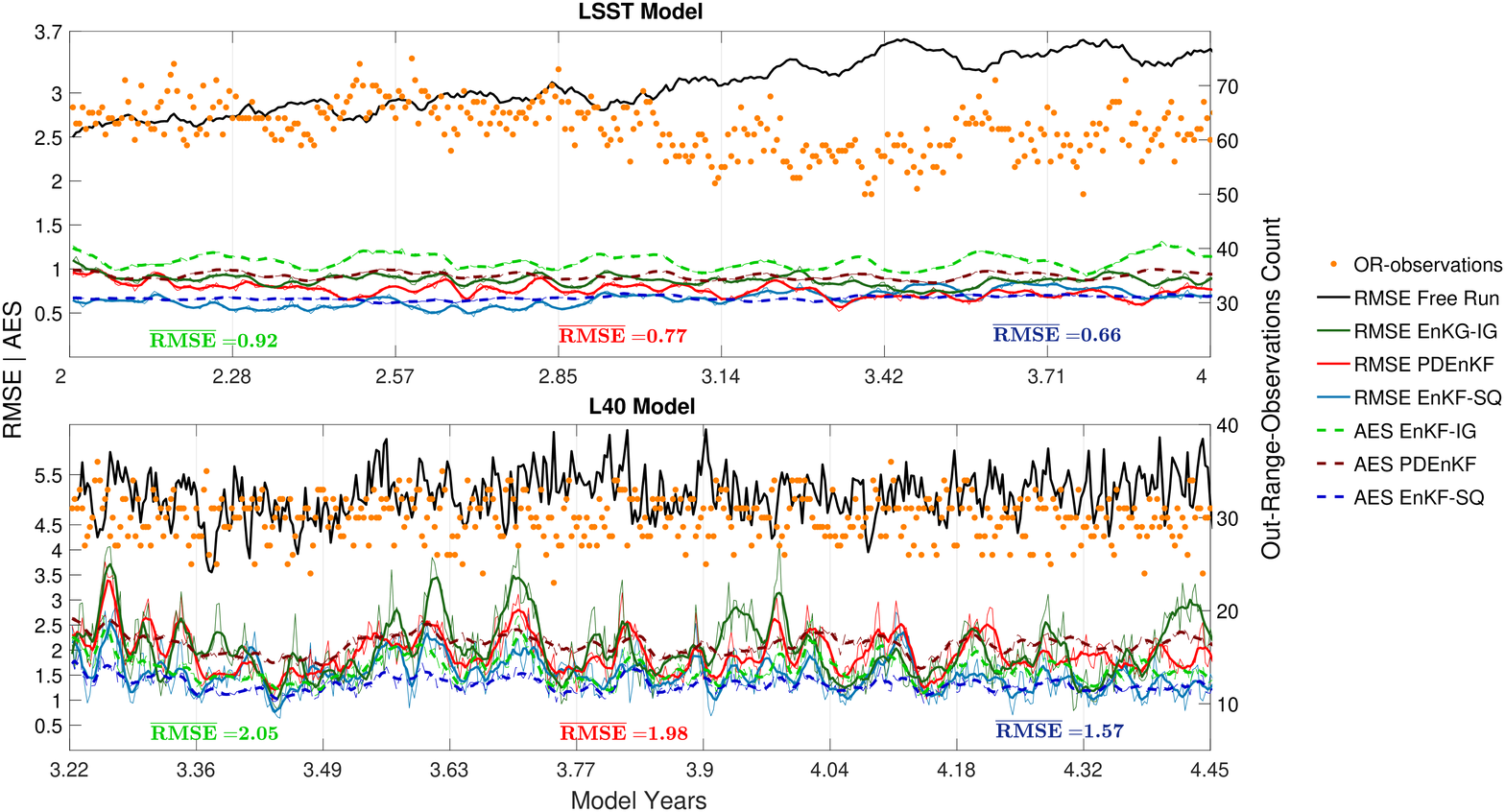}
\caption{Time evolution of the of RMSE (solid lines), AES (dashed lines) and there moving averages in thicker lines. The orange dots represents the number of the OR-observations during the assimilation time. Top and bottom panels show EnKF-SQ, PDEnKF and EnKF-IG results from the LSST and L40 models, respectively. \label{fig:Time-evolution-of-rmse-aes}}
\end{centering}
\end{figure}

Fig.~\ref{fig:Time-evolution-of-rmse-aes} shows the time evolution of RMSE and average ensemble spread (AES) of forecast ensemble members obtained using the EnKF-SQ, PDEnKF and EnKF-IG. Generally, for a ``healthy'' assimilation framework the RMSE is expected to match the AES plus the observation errors. We set different detection limits on observation
in both models such that on average $80\%$ of observations fall out of range i.e., they become soft data. We also show RMSE of a free run (no DA) in both models. For clarity, we superimpose the moving average of RMSE and AES of all $3$ schemes in both panels of the Fig.~\ref{fig:Time-evolution-of-rmse-aes}. As shown in the Fig.~\ref{fig:Time-evolution-of-rmse-aes}, assimilating soft data using the EnKF-SQ improves the forecast RMSE in both models. As shown, among the $3$ tested filters the EnKF-IG is the least accurate. Clearly, assimilating less data degrades the quality of the forecast. We note that the RMSE and total spread (AES + observation error standard deviation) are of the same order, indicating no signs of inbreeding or divergence. As shown, both the EnKF-SQ and the PDEnKF benefit from assimilating soft data. On average, the proposed EnKF-SQ estimates are $20$\% and $12$\% more accurate than those of the PDEnKF for the LSST and L40 models, respectively (Fig.~\ref{fig:Time-evolution-of-rmse-aes}).

\begin{figure}
\begin{centering}
\includegraphics[width=17cm,height=14cm]{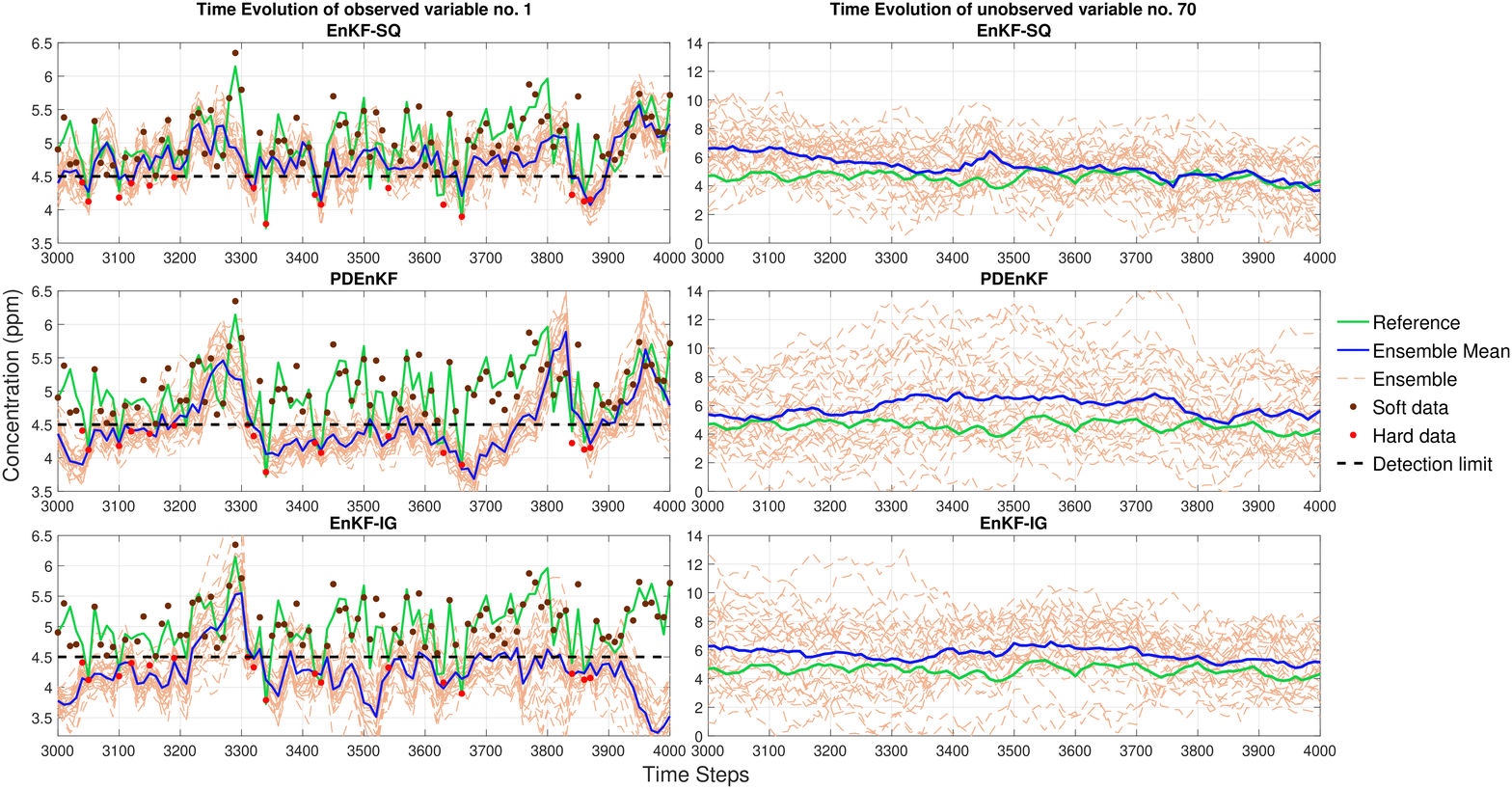}
\caption{Time evolution of the forecast ensemble (dashed orange), ensemble mean (solid blue) and truth (solid green) for observed and unobserved state variable number 1 (left panels) and 70 (right panels), respectively in the LSST model, obtained using EnKF-SQ (top panels), PDEnKF (middle) and EnKF-IG (bottom). \label{fig:Time_evolution_state_variable_lss_model}}
\end{centering}
\end{figure}

To visualize the time evolution of the ensemble for the EnKF-SQ, PDEnKF and EnKF-IG with the LSST model, we plot the concentration of a randomly chosen observed and unobserved state variable versus time in Fig.~\ref{fig:Time_evolution_state_variable_lss_model}. Analyzing the results from Fig.~\ref{fig:Time_evolution_state_variable_lss_model} along with Fig.~\ref{fig:Time-evolution-of-rmse-aes} clearly demonstrates that assimilating soft data not only improve the RMSE, but also reduces the uncertainty in the forecast by shrinking the ensemble spread around the truth. As expected, the EnKF-IG is the least accurate, generating low concentration values when the observations are above the threshold, which the EnKF-SQ avoids successfully. Compared to the PDEnKF, the proposed scheme matches better the truth trajectory. This can be clearly observed for the time intervals (3000, 3700). Similar behavior was also observed for the L40 estimates (not shown).
\subsubsection{Sensitivity Experiments\label{subsec:Sensitivity_experiments} }

\begin{figure}
\subfloat[Sensitivity with variable ensemble size\label{fig:RMSE_L40_different_ens_size}]{\includegraphics[width=8.3cm,height=6.5cm]{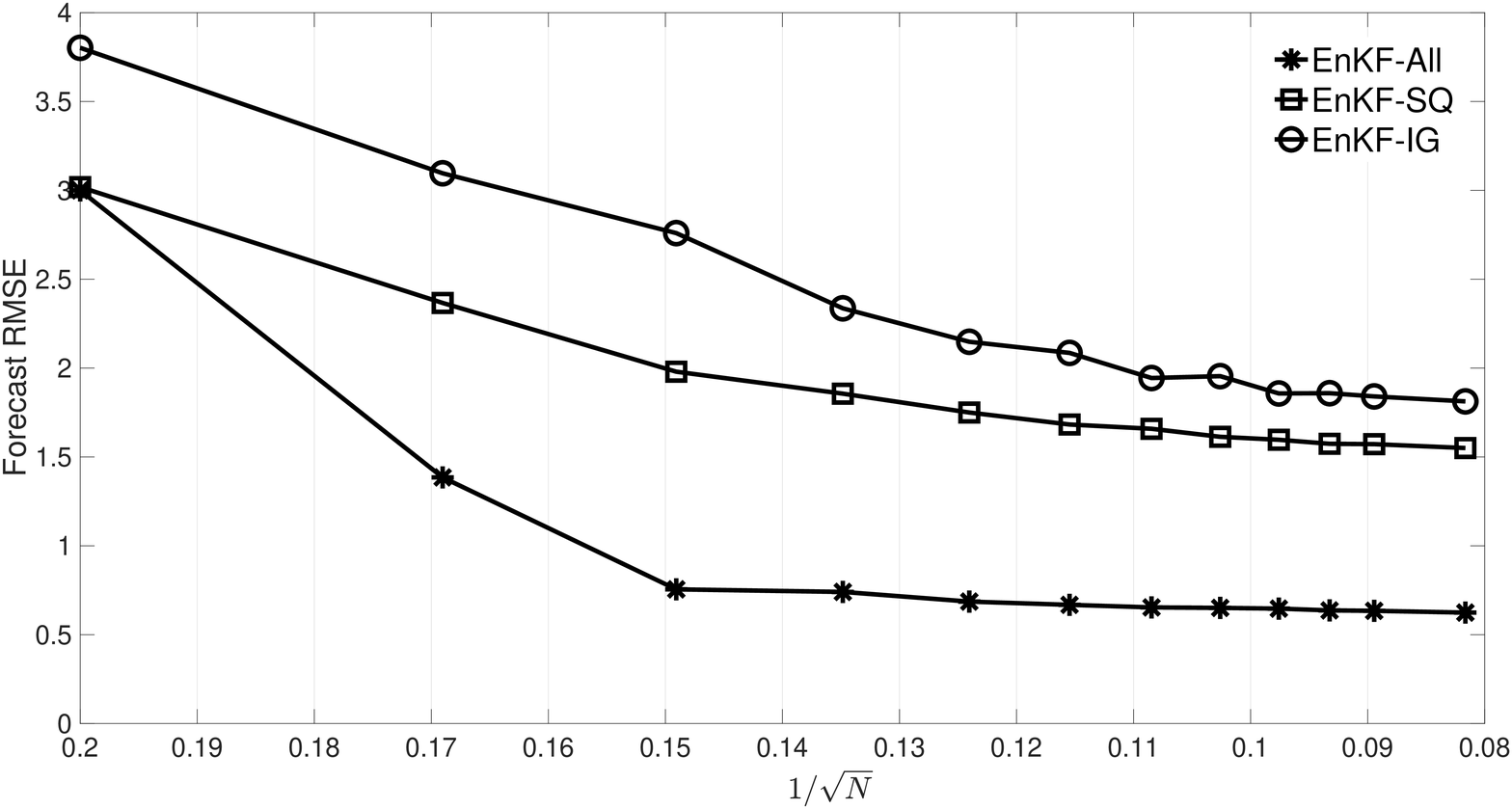} }
\hspace{1em}
\subfloat[ Sensitivity with variable obs. detection limit\label{fig:RMSE_L40_diff_detection_limit}]{\includegraphics[width=8.3cm,height=6.5cm]{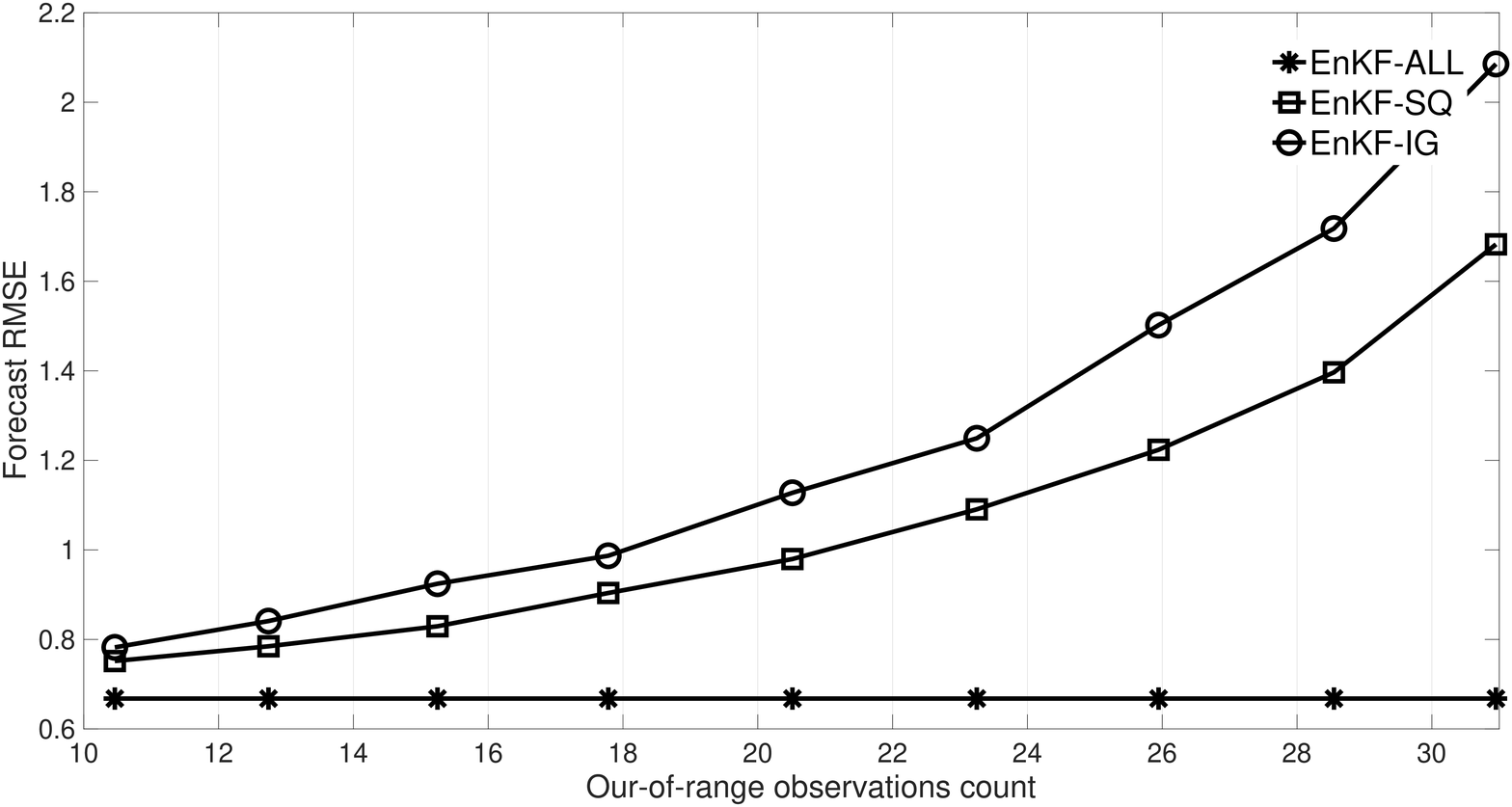}}
\caption{Root mean square error of the forecast estimates resulting from EnKF-All, EnKF-SQ and EnKF-IG for sensitivity experiments with varying ensemble size ($N=25,35,45,...,150$) in left panel and observation detection limit in right panel for L40 model.\label{fig:Sensitivity_experiments_results}}
\end{figure}

Sensitivity experiments are conducted by varying both the ensemble size and detection limits using the L40 model and the results are presented in Fig.~\ref{fig:Sensitivity_experiments_results}. The RMSE values obtained for these experiments are averaged over 5 year-long DA run. The goal is to assess the convergence rate of the EnKF-SQ while increasing $N$ from 25 to 150. The resulting RMSE is plotted against $1/\sqrt{N}$ given that the precision of Monte Carlo methods is varies as a function of $1/\sqrt{N}$. As the ensemble size increases, the RMSE values for each scheme naturally decrease as shown in Fig.~\ref{fig:RMSE_L40_different_ens_size}, although none of them is linear in $1/\sqrt{N}$. For all tested ensemble sizes, the proposed scheme is consistently more accurate than the EnKF-IG. We also note that for small ensemble sizes, the EnKF-SQ performance is as good as the EnKF-ALL and for $N=25$, the estimates of both schemes overlap. 

Changing the detection limit on observations is done such that the number of observations falling out-of-range increases gradually and system has less hard data to assimilate. The forecast RMSE resulting from the EnKF-SQ is shown in (Fig.~\ref{fig:RMSE_L40_diff_detection_limit}) to vary between two extreme cases; i.e., EnKF-IG and EnKF-ALL. Even with a very few hard data to assimilate, EnKF-SQ estimates are alomst $19\%$ more accurate than those of the EnKF-IG. All three schemes converge towards the same RMSE as more hard data is assimilated.
By comparison to the EnKF, $\sigma_{or}$ is the only new parameter introduced in the EnKF-SQ. This imposes only minor changes to existing EnKF codes. We perform sensitivity experiments by introducing a scalar multiplier to eq.~(\ref{eq:sigma_OR}), namely $\alpha$, to examine the behavior of the EnKF-SQ and the impact of using more skewed ensembles. We change $\alpha$ between $0.05$ to $1.85$ with a step size of $0.15$. The new form of eq.~(\ref{eq:sigma_OR}) is shown below:
\begin{equation}
\sigma_{or^{*}}=\alpha\underset{\sigma_{or}}{\underbrace{\left[-\mu+\left(\int\displaylimits_{\mu}^{+\infty}\mathbf{y}f_{\text{clim}}(\mathbf{y})dy\right)\right]}}\label{eq:sensitivity_sig_or}
\end{equation}

\begin{figure}
\subfloat[L40 Model\label{fig:l96_sensitivity_sig_or}]{\includegraphics[width=8.3cm,height=6.5cm]{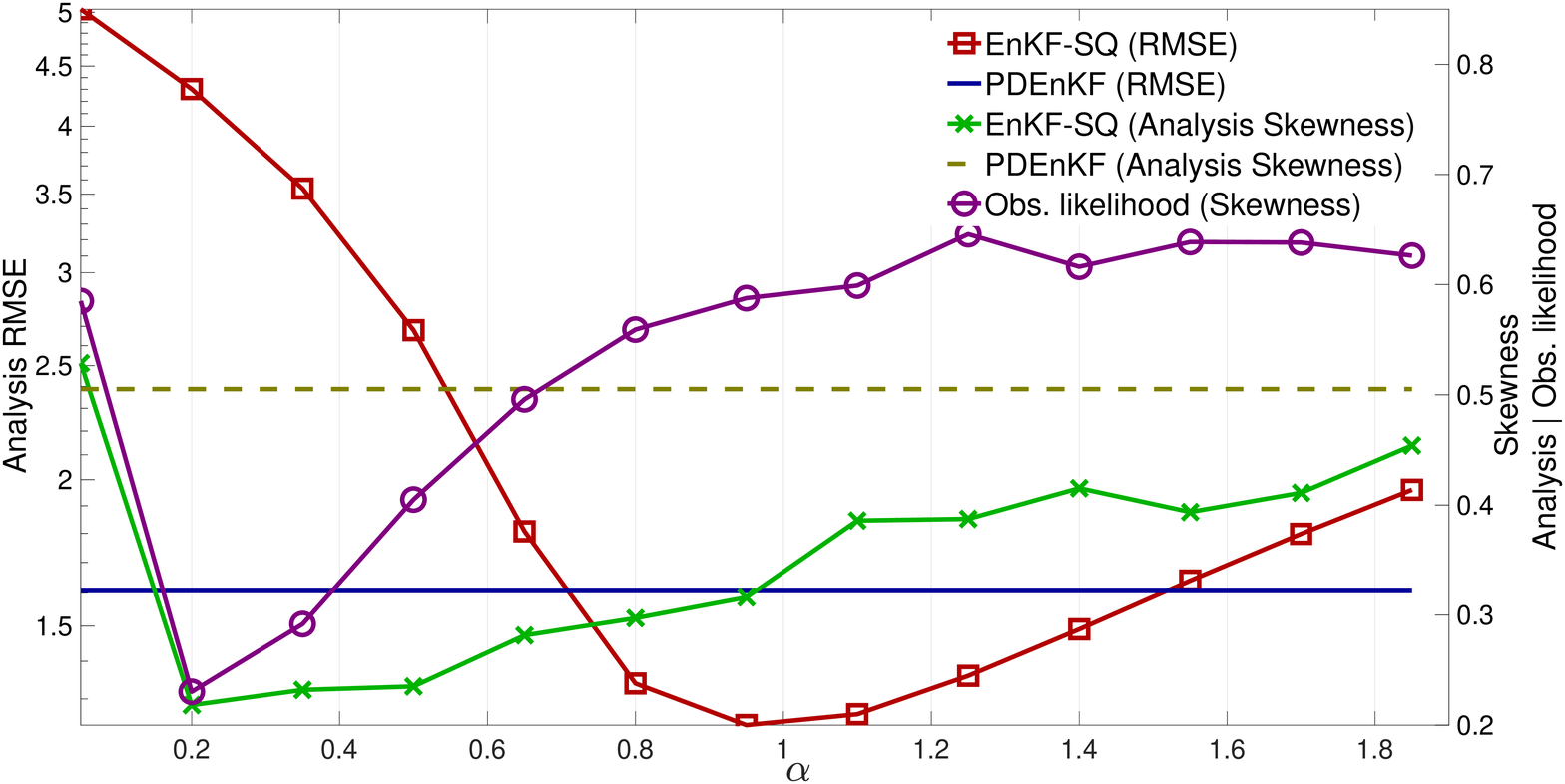}}
\hspace{1em}
\subfloat[LSST Model\label{fig:LSST_sensitivity_sig_or}]{ \includegraphics[width=8.3cm,height=6.5cm]{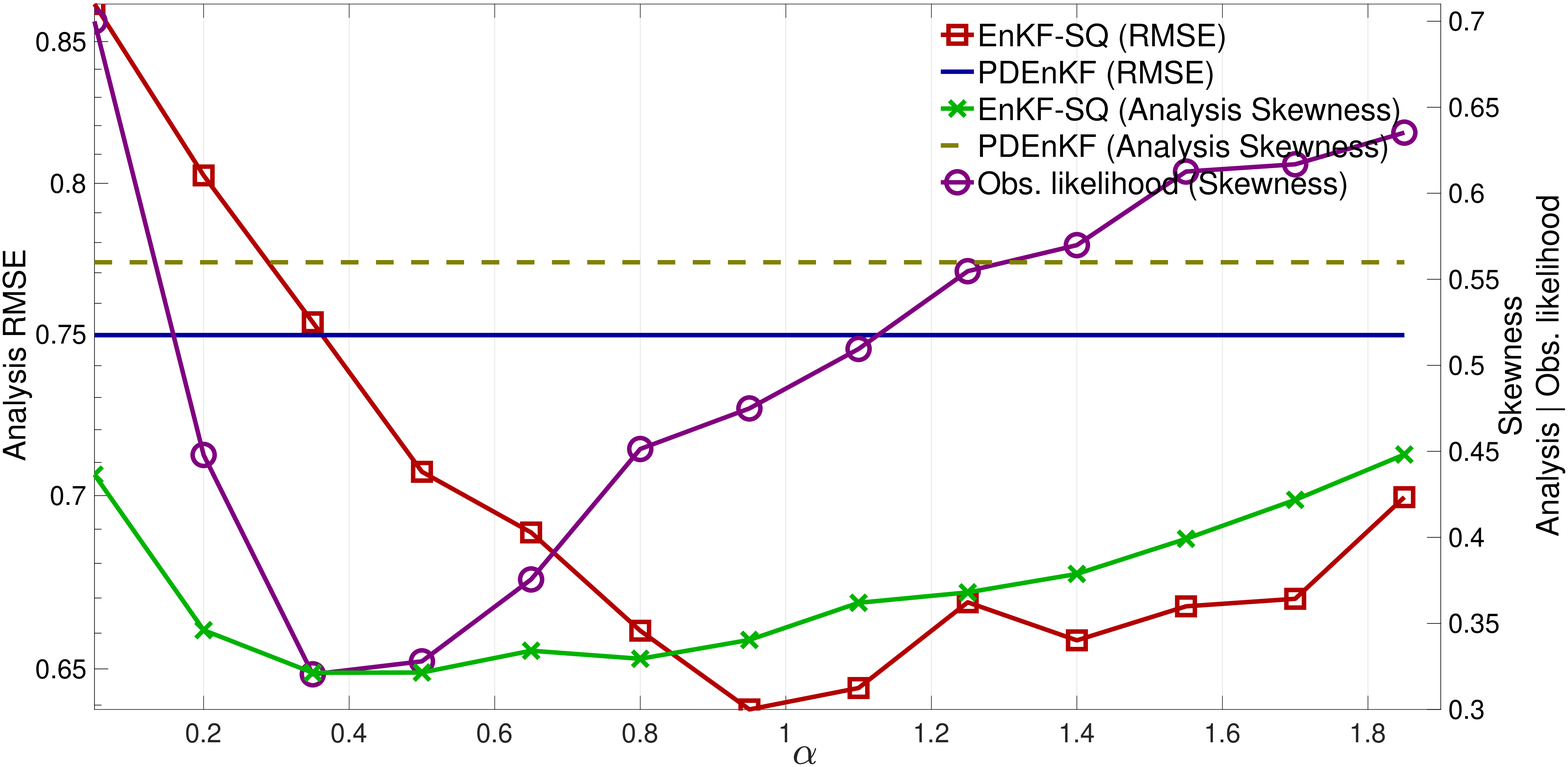}}
\caption{Performance of the EnKF-SQ and PDEnKF using 30 and 75 ensemble members for the LSST (right panel) and L40 (left panel) with increasing values of $\alpha$. The analysis RMSE in addition to the skewness of the posterior and likelihood distributions are demonstrated.\label{fig:Sensitivity_sig_or_both} }
\end{figure}

We plot the RMSE values of the analysis states, from both models, in addition to the absolute skewness of the analysis and observation likelihood versus $\alpha$ for the EnKF-SQ in the Fig.~\ref{fig:Sensitivity_sig_or_both}. The values obtained by the PDEnKF are independent of $\alpha$ and shown for reference. The skewness of analysis and observation likelihood are evaluated as the average absolute value of each variable skewness and only at the last assimilation step: 
\begin{equation}
\text{skew}_{a}=\frac{1}{n}\sum\limits_{j=1}^n\Biggl|\frac{\frac{1}{N}\sum\limits_{i=1}^N\left(\mathbf{x}_{t_{max},j,i}^{a}-\hat{\mathbf{x}}_{t_{max},j}^{a}\right)^{3}}{\left(\frac{1}{N}\sum\limits_{i=1}^N\left(\mathbf{x}_{t_{max},j,i}^{a}-\hat{\mathbf{x}}_{t_{max},j}^{a}\right)^{2}\right)^{3/2}}\Biggl|
\end{equation}

\begin{equation}
\text{skew}_{o}=\frac{1}{m}\sum\limits_{j=1}^m\Biggl|\frac{\frac{1}{N}\sum\limits_{i=1}^N\left(\mathbf{y}_{t_{max},j,i}-\hat{\mathbf{y}}_{t_{max},j}\right)^{3}}{\left(\frac{1}{N}\sum\limits_{i=1}^N\left(\mathbf{y}_{t_{max},j,i}-\hat{\mathbf{y}}_{t_{max},j}\right)^{2}\right)^{3/2}}\Biggl|
\end{equation} 
where $m$ is the number of observation; $\mathbf{x}_{t_{max},j,i}^{a}$ and $\mathbf{y}_{t_{max},j,i}$ are the $i^{th}$ analysis ensemble member and $i^{th}$ observation perturbation vector at time $t_{max}$, respectively. $\hat{\mathbf{x}}_{t_{max},j,i}^{a}$ and $\hat{\mathbf{y}}_{t_{max},j,i}$ are the analysis mean and mean of observation perturbation vector at time $t_{max}$, respectively.

RMSE changes in Fig.~\ref{fig:Sensitivity_sig_or_both}, indicate that when the value of $\alpha$ approaches to 1 i.e., close to nominal value of $\sigma_{or}$ in eq.~(\ref{eq:sigma_OR}), the EnKF-SQ outperforms the PDEnKF. As the value of $\alpha$ moves away from $1$, the performance of EnKF-SQ start to deteriorate, especially for L40 model. This, however, is less obvious for the LSST model.. This can be explained by a poor sampling of the two-piece Gaussian likelihood using a finite ensemble size when $\sigma_{or}$ is assigned very high and/or low values. For instance in the case of a high $\sigma_{or}$ sampling might produce very large perturbations of observations (outliers), which can make the analysis increments more erratic. On the other hand, small values of $\sigma_{or}$ are also detrimental as they are prone to generate samples concentrated around the detection limit, thus pulling the analysis close to an artificial threshold limit. This confirms a posteriori the choice of the nominal value of $\sigma_{or}$ and the importance of a good knowledge of climatological values: if the climatological average of L40 OR-values is biased by more than $50$\% ($\alpha$ lower than $0.5$ or larger than $1.5$), then the flat likelihood of the PDEnKF makes a better option. The linear LSST model is more permissive in this respect since the EnKF-SQ will beat the PDEnKF even with values of $\alpha$ more than $100$\% off the nominal value. This could be because non-Gaussianity is reduced in a linear model like the LSST (central limit theorem), contrary to the nonlinear and chaotic L40 model. 

Fig.~\ref{fig:Sensitivity_sig_or_both} also shows the absolute skewness of the observation likelihood and analysis ensemble. The skewness of the EnKF-SQ analysis ensemble follows the same trend as that of the likelihood, though the linear EnKF-SQ update makes it less skewed. To illustrate, as $\alpha$ increases, the observation likelihood transitions from being right to left skewed (not shown). Likewise the analysis follows a similar behavior. The analysis ensembles are quite severely skewed (typical skewness is from $0.3$ to $0.5$ for the EnKF-SQ and higher with the PDEnKF update scheme), which does not seem to affect the EnKF-SQ performance directly. The minimum RMSE does not even coincide with the minimum skewness. This indicates that the method can handle some degree of non-Gaussianity, which makes it useful for assimilating soft data with the EnKF-SQ and PDEnKF. 

\section{Summary and discussion\label{sec:Conclusion}}

Many observations in practice are only available within a confined range. Qualitative information measured above or below the detection limit can still be exploited by data assimilation although current methods only consider hard data. In this paper, we proposed a new DA algorithm referred to as EnKF-SQ in order to assimilate semi-qualitative observations through an explicit treatment of soft data. The update algorithm requires a pre-process step, in which observations are split into two groups, hard and soft data. This is then followed by the update of forecast ensemble using the Kalman update. An assumption is imposed on the observation likelihood to be a two-piece Gaussian and the mode of the likelihood is positioned at the detection limit. Members falling inside or outside the observable range are then separated to have consistent update by soft data. This makes it necessary to update each forecast ensemble member individually, but not in parallel. Computationally, this is not a major turn-down as in many applications the update only represents a few percents of the costs of the ensemble propagation step \citep{Sakov2012} and a local EnKF-SQ would still run local updates in a parallel loop. 

The new EnKF-SQ has been evaluated in a linear subsurface transport and a nonlinear Lorenz-40 model. Its performance has been compared two different versions of the stochastic EnKF, namely EnKF-ALL (no detection limit on observations) and EnKF-IG (no assimilation of soft data) in addition to the previously introduced partial deterministic ensemble Kalman filter (PDEnKF), which is built over a deterministic EnKF and uses a uniform OR prior likelihood. Our numerical results suggest that assimilating soft data with the EnKF-SQ improves the overall forecast accuracy. The scheme outperforms the EnKF-IG with reasonable computing time and ensemble sizes lower than $100$ for systems of dimension greater than $20$. Thus it does
not suffer from the curse of dimensionality. This suggests that EnKF-SQ is a viable method and can be implemented with more realistic applications of the EnKF.

Sensitivity experiment to the chosen value of OR-observation likelihood error variance $\sigma_{or}$ imply that, if chosen \textit{properly}, the EnKF-SQ performs better than the PDEnKF. This may not be true for all types of applications however, because the differences in performance are small and some observations may be represented by an OR-likelihood with a fatter tail than the Gaussian distribution. Such cases are not addressed in the present work. As far as the two-piece Gaussian likelihood goes, we found that even though it might increase the skewness of the posterior distribution, the benefits of assimilating soft data out-beat the inconvenience of non-Gaussianity in both linear and non-linear cases. 

The question arises whether the assimilation of an arbitrary value in the out-of-range domain will perform as good as the EnKF-SQ, with less algorithmic complexity. This has not been tested but we note that assimilating a hard pseudo-observation in the OR domain would not introduce asymmetric information as the EnKF-SQ does, so the approach
would unnecessarily update forecast ensemble members that fall rightly in the OR domain. 

Is this semi-qualitative approach applicable to other data assimilation methods? It requires a stochastic data assimilation method to treat the ensemble members as possible realizations of the underlying random variables. Extensions to deterministic methods are therefore not straightforward, the link to Optimal Interpolation (OI) can be done by geostatistical methods through randomization \citep{Emery2008}, but this would make the OI method much more costly. The extension from ensemble filters to ensemble smoothers should however be straightforward. 

\section{Acknowledgements}

The Authors would like to thank Morten Borup for interesting discussions and hosting AS at DTU. We would also like to thank Fran\c{c}ois Counillon for the insightful scientific discussions and Alberto Carrassi for suggesting the name `EnKF-SQ'. The research is funded by the Nordic Center of Excellence EmblA (Ensemble-based data assimilation for environmental monitoring and prediction) under NordForsk contract number 56801.

\bibliographystyle{wileyqj}
\bibliography{EnKF_SQ_references}

\end{document}